\newcount\ite\ite=1\def\0{\global\ite=1\1}
\def\1{\item{\rm(\romannumeral\the\ite)}\advance\ite1\quad}
\def\phi{\varphi}

\documentclass[12pt,twoside]{article}
\usepackage{amssymb}

\font\teneufm=eufm10 scaled \magstep1
\font\seveneufm=eufm7 scaled \magstep1
\font\fiveeufm=eufm5  scaled \magstep1
\newfam\eufmfam

\textfont\eufmfam=\teneufm
\scriptfont\eufmfam=\seveneufm
\scriptscriptfont\eufmfam=\fiveeufm
\def\frak#1{{\fam\eufmfam\relax#1}}

\newfam\msbfam
\font\tenmsb=msbm10 scaled \magstep1  \textfont\msbfam=\tenmsb
\font\sevenmsb=msbm7 scaled \magstep1 \scriptfont\msbfam=\sevenmsb
\font\fivemsb=msbm5 scaled \magstep1  \scriptscriptfont\msbfam=\fivemsb

\makeatletter
\def\blfootnote{\xdef\@thefnmark{}\@footnotetext}
\makeatother

\def\dd#1{\raise1.5pt\hbox{$\,\partial\!$}/\raise-2.5pt\hbox{$\!\partial#1\,$}}

\def\tilde{\widetilde}

\def\5#1{{\mathcal #1}}

\def\CC{{\mathbb C}}

\def\ra{\rightarrow}

\def\GL{\mathop{\rm GL}\nolimits}

\def\Ann{\mathop{\rm Ann}\nolimits}

\def\emb{\mathop{\rm emb}\nolimits}
\def\Soc{\mathop{\rm Soc}\nolimits}
\def\char{\mathop{\rm char}\nolimits}

 \def\HollowBoxx #1#2#3{{\dimen0=#1 \advance\dimen0 by -#2
       \dimen1=#1 \advance\dimen1 by #3
        \vrule height 0pt depth #3 width #2
       \hskip -#3
       \vrule height #1 depth #3 width #3}}
 \def\LeftContraction{\mathord{\kern1.45pt \HollowBoxx{6pt}{3.5pt}{.4pt}}\,}

 \def\HollowBox #1#2#3{{\dimen0=#1 \advance\dimen0 by -#3
       \dimen1=#1 \advance\dimen1 by #3
        \vrule height #1 depth #3 width #3
        \vrule height 0pt depth #3 width #2
        \hskip -#3}}
 \def\RightContraction{\mathord{\, \HollowBox{6pt}{3.1pt}{.4pt}} \kern1.6pt}

\def\qed{{\hfill $\Box$}}
\newtheorem{theorem}{THEOREM}[section]
\newtheorem{corollary}[theorem]{Corollary}
\newtheorem{lemma}[theorem]{Lemma}
\newtheorem{example}[theorem]{Example}
\newtheorem{remark}[theorem]{Remark}

\begin{document}

\begin{center}
{\Large \bf A Criterion for Isomorphism\\
\vspace{0.3cm}
of Artinian Gorenstein Algebras}\blfootnote{{\bf Mathematics Subject Classification:} 13H10}\blfootnote{{\bf Keywords:} Artinian Gorenstein algebras}
\medskip\\
\normalsize A. V. Isaev
\end{center}

\begin{quotation} 
{\small \sl \noindent 
Let $A$ be an Artinian Gorenstein algebra over an infinite field $k$ with either $\char(k)=0$ or $\char(k)>\nu$, where $\nu$ is the socle degree of $A$. To every such algebra and a linear projection $\pi$ on its maximal ideal ${\mathfrak m}$ with range equal to the socle $\Soc(A)$ of $A$, one can associate a certain algebraic hypersurface $S_{\pi}\subset{\mathfrak m}$, which is the graph of a polynomial map $P_{\pi}:\ker\pi\ra\Soc(A)\simeq k$. Recently, in {\rm\cite{FIKK}}, {\rm\cite{FK}} the following surprising criterion was obtained: two Artinian Gorenstein algebras $A$, $\tilde A$ are isomorphic if and only if any two hypersurfaces $S_{\pi}$ and $S_{\tilde\pi}$ arising from $A$ and $\tilde A$, respectively, are affinely equivalent. The proof is indirect and relies on a geometric argument. In the present paper we give a short algebraic proof of this statement. We also discuss a connection, established elsewhere, between the polynomials $P_{\pi}$ and Macaulay inverse systems.}
\end{quotation}

\thispagestyle{empty}

\pagestyle{myheadings}
\markboth{A. V. Isaev}{Criterion for Isomorphism of Artinian Gorenstein Algebras}

\setcounter{section}{0}

\section{Introduction}\label{intro}
\setcounter{equation}{0}

We consider Artinian local commutative associative algebras over a field $k$. Recall that such an algebra $A$ is Gorenstein if and only if the socle $\Soc(A)$ of $A$ is a 1-dimensional vector space over $k$ (see, e.g.~\cite{Hu}). Gorenstein algebras frequently occur in various areas of mathematics and its applications to physics (see, e.g.~\cite{B}, \cite{L}). For $k=\CC$, in \cite{FIKK} we found a surprising criterion for two Artinian Gorenstein algebras to be isomorphic. The criterion was given in terms of a certain algebraic hypersurface $S_{\pi}$ in the maximal ideal ${\mathfrak m}$ of $A$ associated to a linear projection $\pi$ on ${\mathfrak m}$ with range $\Soc(A)$, where we assume that $\dim_k A>1$. The hypersurface\linebreak $S_{\pi}$ passes through the origin and is the graph of a polynomial map\linebreak $P_{\pi}:\ker\pi\ra \Soc(A)\simeq\CC$. In \cite{FIKK} we showed that for $k=\CC$ two Artinian Gorenstein algebras $A$, $\tilde A$ are isomorphic if and only if any two hypersurfaces $S_{\pi}$ and $S_{\tilde\pi}$ arising from $A$ and $\tilde A$, respectively, are affinely equivalent, that is, there exists a bijective affine map ${\mathcal A}: {\mathfrak m}\ra\tilde {\mathfrak m}$ such that ${\mathcal A}(S_{\pi})=S_{\tilde\pi}$.

It should be noted that the above criterion differs from the well-known criterion in terms of Macaulay inverse systems, where to an Artinian Gorenstein algebra $A$ one also associates certain polynomials (see, e.g.~Proposition 2.2 in \cite{ER} and references therein). Indeed, as the discussion in Section \ref{invsys} shows, in general none of the polynomials $P_{\pi}$ is an inverse system for $A$. The value of the method of \cite{FIKK} lies in the fact that affine equivalence for the hypersurfaces $S_{\pi}$ and $S_{\tilde\pi}$ is sometimes easier to verify than the kind of equivalence required for inverse systems (see formula (7) in \cite{ER}). In Section \ref{app} we give an example when this is indeed the case.

At the time of writing the article \cite{FIKK} we could not find an algebraic proof of our criterion, but, quite unexpectedly, discovered a proof based on the geometry of certain real codimension two quadrics in complex space. On the other hand, the hypersurface $S_{\pi}$ can be introduced for an Artinian Gorenstein algebra $A$ over any field $k$ of characteristic either zero or greater than the socle degree $\nu$ of $A$. Therefore, it is natural to attempt to extend the result of \cite{FIKK} to all such algebras. Indeed, an extension of this kind was obtained in the recent paper \cite{FK}. The result was stated for fields of zero characteristic, but the proof is likely to work at least for some fields satisfying $\char(k)>\nu$ as well. The argument proceeds by reducing the case of an arbitrary field to the case $k=\CC$, and  therefore the basis of the method of \cite{FK} remains the geometric idea of \cite{FIKK}. It would be very desirable, however, to find a purely algebraic argument directly applicable to any field. In this paper we present such an argument, in which neither geometry nor reduction to the case $k=\CC$ is required and which works for any infinite field $k$ with either $\char(k)=0$ or $\char(k)>\nu$. We note that a very brief outline of the argument was given in \cite{Is}, and the main purpose of this article is to provide full details. Another purpose of the paper is to exhibit our approach to Artinian Gorenstein algebras to a broad audience, and therefore the article is to some extent expository.

The paper is organized as follows. Section \ref{prelim} contains necessary preliminaries and the precise statement of the criterion in Theorem \ref{equivalence}. Our proof of Theorem \ref{equivalence} is given in Section \ref{proofs1}. In Section \ref{app} we demonstrate how this result works in applications (see Example \ref{appl}). Namely, we consider an Artinian Gorenstein algebra $A_0$, originally presented in \cite{FK}, of embedding dimension 3 and socle degree 4 whose Hilbert function is $\{1,3,3,1,1\}$. Further, we perturb $A_0$ to obtain a one-parameter family $A_t$ of Artinian Gorenstein algebras of the same dimension and socle degree. While directly establishing a relationship between $A_t$ and $A_0$ seems to be nontrivial, this can be easily done with the help of Theorem \ref{equivalence}. Namely, using the theorem we show, in a rather elementary way, that each algebra $A_t$ is in fact isomorphic to $A_0$. 

Next, in Section \ref{invsys} we discuss the connection between the polynomials $P_{\pi}$ and Macaulay inverse systems for any Artinian Gorenstein algebra $A$ found in our earlier article \cite{AI}. The proof is short, and for the completeness of our exposition we include it in the present paper as well. Namely, if $P_{\pi}$ is regarded as a map from $\ker\pi$ to $k$ (in which case we call it a {\it nil-polynomial}), then the restriction of $P_{\pi}$ to any subspace of $\ker\pi$ that forms a complement to ${\mathfrak m}^2$ in ${\mathfrak m}$ is an inverse system for $A$ (see Remark \ref{nilpolinvsys}). All these subspaces are of dimension equal to the embedding dimension $\emb\dim A$ of $A$, thus nil-polynomials can be viewed as extensions of certain inverse systems to spaces of dimension greater than $\emb\dim A$. As a result, since nil-polynomials are given explicitly, one obtains an effective formula for computing an inverse system for any Artinian Gorenstein algebra (see Theorem \ref{main} and Remark \ref{remembdim}). It appears that such a formula had not existed in the literature prior to our work.

Utilizing the above relationship between nil-polynomials and Mac\-au\-lay inverse systems, at the end of Section \ref{invsys} we revisit Example \ref{appl} and explain how isomorphism between $A_t$ and $A_0$ can be established by using inverse systems. This turns out to be significantly harder to do than by applying the method based on nil-polynomials. In general, to use inverse systems one needs to compute them first, and Theorem \ref{main} provides an explicit way for producing them from nil-polynomials. This fact alone shows that nil-polynomials are a new rather useful tool for the study of Artinian Gorenstein algebras.   

{\bf Acknowledgements.} We are grateful to N. Kruzhilin for stimulating discussions and to the referee for the thorough reading of the paper and helpful suggestions. This work is supported by the Australian Research Council.

\section{Preliminaries}\label{prelim}
\setcounter{equation}{0}
    
Let $A$ be an Artinian Gorenstein algebra over a field $k$ with identity element ${\mathbf 1}$, maximal ideal ${\mathfrak m}$ and socle $\Soc(A):=\{x\in A: x\,{\mathfrak m}=0\}$. Suppose that $\dim_kA>1$ (i.e.~${\mathfrak m}\ne 0$) and denote by $\nu$ the socle degree of $A$, that is, the largest among all integers $\mu$ for which ${\frak m}^{\mu}\ne 0$. Observe that $\nu\ge 1$ and $\Soc(A)={\mathfrak m}^{\nu}$. 

Assume now that either $\char(k)=0$ or $\char(k)>\nu$ and define the {\it exponential map}\, $\exp: {\mathfrak m}\ra {\bf 1}+{\mathfrak m}$ by the formula
$$
\displaystyle\exp(x):=\sum_{m=0}^{\nu}\frac{1}{m!}x^m,
$$
where $x^0:={\bf 1}$. This map is bijective with the inverse given by
$$
\log({\bf 1}+x):=\sum_{m=1}^{\nu}\frac{(-1)^{m+1}}{m}x^m,\quad x\in{\mathfrak m}.
$$

Next, fix a linear projection $\pi$ on $A$ with range $\Soc(A)$ and kernel containing ${\bf 1}$ (we call such projections {\it admissible}). Set ${\mathcal K}:=\ker\pi\hspace{0.01cm}\cap\hspace{0.01cm}{\mathfrak m}$ and let $S_{\pi}$ be the graph of the polynomial map $P_{\pi}:{\mathcal K}\ra\Soc(A)$ of degree $\nu$ defined as follows:
\begin{equation}
P_{\pi}(x):=\pi(\exp(x))=\pi\left(\sum_{m=2}^{\nu}\frac{1}{m!}x^m\right),\quad x\in{\mathcal K}\label{poly}
\end{equation}
(note that for $\dim_k A=2$ one has $P_{\pi}=0$). Everywhere below we identify any point $(x,P_{\pi}(x))\in S_{\pi}$ with the point $x+P_{\pi}(x)\in{\mathcal K}\oplus\Soc(A)={\mathfrak m}$ and therefore think of $S_{\pi}$ as a hypersurface in ${\mathfrak m}$. Observe that the $\Soc(A)$-valued quadratic part of $P_{\pi}$ is non-degenerate on ${\mathcal K}$ since the $\Soc(A)$-valued bilinear form
\begin{equation}
b_{\pi}(a,c):=\pi(ac), \quad a,c\in A\label{bpi}
\end{equation}
is non-degenerate on $A$ (see, e.g.~p.~11 in \cite{He}). Examples of hypersurfaces $S_{\pi}$ explicitly computed for particular algebras can be found in \cite{FIKK}, \cite{FK}, \cite{EI} (see also Section \ref{app} below).

We will now state the criterion for isomorphism of Gorenstein algebras that was obtained in \cite{FIKK}, \cite{FK} for algebras over fields of zero characteristic.

\begin{theorem}\label{equivalence}\sl Let $A$, $\tilde A$ be Gorenstein algebras of finite vector space dimension greater than 1 and socle degree $\nu$ over an infinite field $k$ with either $\char(k)=0$ or $\char(k)>\nu$. Let, further, $\pi$, $\tilde\pi$ be admissible projections on $A$, $\tilde A$, respectively. Then $A$, $\tilde A$ are isomorphic if and only if the hypersurfaces $S_{\pi}$, $S_{\tilde\pi}$ are affinely equivalent. Moreover, if ${\mathcal A}:{\mathfrak m}\ra\tilde{\mathfrak m}$ is a linear isomorphism such that ${\mathcal A}(S_{\pi})=S_{\tilde\pi}$, then ${\mathcal A}$ is an algebra isomorphism.
\end{theorem}

\section{Proof of Theorem \ref{equivalence}}\label{proofs1}
\setcounter{equation}{0}

For every hypersurface $S_{\pi}$, we let ${\mathcal S}_{\pi}$ be the graph over ${\mathcal K}$ of the polynomial map $-P_{\pi}$ (see (\ref{poly})). Observe that
\begin{equation}
{\mathcal S}_{\pi}=\{x\in{\mathfrak m}:\pi(\exp(x))=0\}.\label{tildes}
\end{equation}
Clearly, $S_{\pi}$ and $S_{\tilde\pi}$ are affinely equivalent if and only if ${\mathcal S}_{\pi}$ and ${\mathcal S}_{\tilde\pi}$ are affinely equivalent, and below we will obtain the first statement of the theorem with ${\mathcal S}_{\pi}$ and ${\mathcal S}_{\tilde\pi}$ in place of $S_{\pi}$ and $S_{\tilde\pi}$, respectively.

The necessity implication is proved as in Proposition 2.2 of \cite{FIKK}. The argument is short, and for the completeness of our exposition we reproduce it below. The idea is to show that if $\pi_1$, $\pi_2$ are admissible projections on $A$, then ${\mathcal S}_{\pi_{{}_1}}={\mathcal S}_{\pi_{{}_2}}+x_0$ for some $x_0\in{\mathfrak m}$. Clearly, the necessity implication is a consequence of this fact.

For every $y\in{\mathfrak m}$, let $M_y$ be the multiplication operator from $A$ to ${\mathfrak m}$
defined by $a\mapsto ya$ and set ${\mathcal K}_1:=\ker\pi_1\hspace{0.01cm}\cap\hspace{0.01cm}{\mathfrak m}$. The correspondence
$$
y\mapsto
\pi_1\circ M_y|_{{\mathcal K}_1}
$$
defines a linear map ${\mathfrak L}$ from ${\mathcal K}_1$ into the space $L({\mathcal K}_1,\Soc(A))$ of linear maps from ${\mathcal K}_1$ to $\Soc(A)$. Since for every admissible projection $\pi$ the form $b_{\pi}$ defined in (\ref{bpi}) is non-degenerate on $A$ and since $\dim_kL({\mathcal K}_1,\Soc(A))=\dim_k{\mathcal K}_1$, it follows that ${\mathfrak L}$ is an isomorphism. 

Next, let $\lambda:=\pi_2-\pi_1$ and observe that $\lambda({\bf 1})=0$, $\lambda(\Soc(A))=0$. Clearly, $\lambda|_{{\mathcal K}_1}$ lies in $L({\mathcal K}_1,\Soc(A))$, and therefore there exists $y_0\in{\mathcal K}_1$ such that $\lambda|_{{\mathcal K}_1}=\pi_1\circ M_{y_0}|_{{\mathcal K}_1}$. We then have $\lambda=\pi_1\circ M_{y_0}$ everywhere on $A$, hence
$$
\pi_2(\exp(x))=\pi_1\Bigl(({\bf 1}+y_0)\exp(x))\Bigr)=\pi_1(\exp(x+x_0))
$$ 
for $x_0:=\log({\bf 1}+y_0)$, which implies ${\mathcal S}_{\pi_{{}_1}}={\mathcal S}_{\pi_{{}_2}}+x_0$ as claimed.

We will now obtain the sufficiency implication. Let ${\mathcal A}:{\mathfrak m}\ra\tilde{\mathfrak m}$ be an affine equivalence with ${\mathcal A}({\mathcal S}_{\pi})={\mathcal S}_{\tilde\pi}$, and $z_0:={\mathcal A}(0)$. Consider the linear map ${\mathcal L}(x):={\mathcal A}(x)-z_0$, with $x\in{\mathfrak m}$. We will show that ${\mathcal L}:{\mathfrak m}\ra\tilde{\mathfrak m}$ is an algebra isomorphism, which will imply that $A$ and $\tilde A$ are isomorphic. 

Clearly, ${\mathcal L}$ maps ${\mathcal S}_{\pi}$ onto ${\mathcal S}_{\tilde\pi}-z_0$. Consider the admissible projection on $\tilde A$ given by the formula $\tilde\pi'(a):=\tilde\pi(\widetilde{\exp}(z_0)a)$, where $\widetilde{\exp}$ is the exponential map associated to $\tilde A$. Formula (\ref{tildes}) then implies ${\mathcal S}_{\tilde\pi}-z_0={\mathcal S}_{\tilde\pi'}$, hence ${\mathcal L}$ maps ${\mathcal S}_{\pi}$ onto ${\mathcal S}_{\tilde\pi'}$. 

Recall that ${\mathcal S}_{\pi}$ is the graph of the polynomial map $-P_{\pi}:{\mathcal K}\ra\Soc(A)$ (see (\ref{poly})). Set $n:=\dim_k{\mathfrak m}-1=\dim_k\tilde{\mathfrak m}-1$ and choose coordinates\linebreak $\alpha=(\alpha_1,\dots,\alpha_n)$ in ${\mathcal K}$ and a coordinate $\alpha_{n+1}$ in $\Soc(A)$. In these coordinates the hypersurface ${\mathcal S}_{\pi}$ is written as
$$
\alpha_{n+1}=\sum_{i,j=1}^ng_{ij}\alpha_i\alpha_j+\sum_{i,j,\ell=1}^nh_{ijk}\alpha_i\alpha_j\alpha_{\ell}+\cdots,
$$ 
where $g_{ij}$ and $h_{ij\ell}$ are symmetric in all indices, $(g_{ij})$ is non-degenerate, and the dots denote the higher-order terms. In Proposition 2.10 in \cite{FIKK} (which works over any field of characteristic either zero or greater than $\nu$) we showed that the above equation of ${\mathcal S}_{\pi}$ is in {\it Blaschke normal form}, that is, one has $\sum_{ij=1}^ng^{ij}h_{ij\ell}=0$ for all $\ell$, where $(g^{ij}):=(g_{ij})^{-1}$.

We will now need the following lemma, which is a variant of the second statement of Proposition 1 in \cite{EE}.

\begin{lemma}\label{Blashke} \sl Let $V$, $W$ be vector spaces over an infinite field $k$, with\linebreak $\dim_kV=\dim_kW=N+1$, $N\ge 0$. Choose coordinates $\beta=(\beta_1,\dots,\beta_N)$, $\beta_{N+1}$ in $V$ and coordinates $\gamma=(\gamma_1,\dots,\gamma_N)$, $\gamma_{N+1}$ in $W$. Let $S\subset V$, $T\subset W$ be hypersurfaces given, respectively, by the equations
\begin{equation}
\beta_{N+1}={\mathcal P}(\beta_1,\dots,\beta_N),\quad \gamma_{N+1}={\mathcal Q}(\gamma_1,\dots,\gamma_N),\label{sim}
\end{equation}
where ${\mathcal P}$, ${\mathcal Q}$ are polynomials without constant and linear terms. Assume further that equations {\rm (\ref{sim})} are in Blaschke normal form. Then every bijective linear transformation of $V$ onto $W$ that maps $S$ into $T$ has the form
$$
\gamma=C\beta,\quad \gamma_{N+1}= c\,\beta_{N+1},
$$
where $C\in\GL(N,k)$, $c\in k^*$ and $\beta,\gamma$ are viewed as column-vectors.
\end{lemma}

\noindent{\bf Proof:} Let $L:V\ra W$ be a bijective linear transformation. Write $L$ in the most general form 
$$
\gamma=C\beta+d\beta_{N+1},\quad \gamma_{N+1}=\sum_{i=1}^N c_i\beta_i+ c\,\beta_{N+1}
$$
for some $c_1,\dots,c_N,c\in k$, $d\in k^N$, and $N\times N$-matrix $C$ with entries in $k$. Then the condition $L(S)\subset T$ is expressed as 
\begin{equation}
\sum_{i=1}^N c_i\beta_i+c{\mathcal P}(\beta)\equiv {\mathcal Q}(C\beta+d{\mathcal P}(\beta)).\label{equal4}
\end{equation} 

Since $k$ is an infinite field, identity (\ref{equal4}) implies that the coefficients at the same monomials on the left and on the right are equal. Then, since ${\mathcal P}$ and ${\mathcal Q}$ do not contain linear terms, it follows that $c_1=\dots=c_N=0$ and therefore $C\in\GL(N,k)$, $c\in k^*$. Further, comparing the second- and third-order terms in (\ref{equal4}), it is straightforward to see that the equations of both $S$ and $T$ cannot be in Blaschke normal form unless $d=0$. \qed
\vspace{-0.1cm}\\

By Lemma \ref{Blashke}, writing the hypersurface ${\mathcal S}_{\tilde\pi'}$ in some coordinates\linebreak $\tilde\alpha=(\tilde\alpha_1,\dots,\tilde\alpha_n)$, $\tilde\alpha_{n+1}$ in $\tilde{\mathfrak m}$ chosen as above, we see that the map  ${\mathcal L}$ has the form
\begin{equation}
\quad \tilde\alpha= C\alpha,\quad \tilde\alpha_{n+1}=c\,\alpha_{n+1}\label{coord}
\end{equation}
for some $C\in\GL(n,k)$, $c\in k^*$. In coordinate-free formulation, (\ref{coord})\linebreak means that with respect to the decompositions ${\mathfrak m}={\mathcal K}\oplus\Soc(A)$ and\linebreak $\tilde{\mathfrak m}=\tilde{\mathcal K}\oplus\Soc(\tilde A)$, where $\tilde{\mathcal K}:=\ker\tilde\pi'\cap\tilde{\mathfrak m}$, the map ${\mathcal L}$ has the    block-\linebreak diagonal form, that is, there exist linear isomorphisms $L_1: {\mathcal K}\ra\tilde{\mathcal K}$ and $L_2:\Soc(A)\ra\Soc(\tilde A)$ such that ${\mathcal L}(x+u)=L_1(x)+L_2(u)$, with $x\in{\mathcal K}$, $u\in\Soc(A)$. Therefore, for the corresponding polynomial maps $P_{\pi}$ and $P_{\tilde\pi'}$ (see (\ref{poly})) we have
\begin{equation}
L_2\circ P_{\pi}=P_{\tilde\pi'}\circ L_1.\label{lineq1}
\end{equation}
Clearly, identity (\ref{lineq1}) yields
\begin{equation}
L_2\circ P_{\pi}^{[m]}=P_{\tilde\pi'}^{[m]}\circ L_1,\quad m=2,\dots,\nu,\label{lineq}
\end{equation}
where $P_{\pi}^{[m]}$, $P_{\tilde\pi'}^{[m]}$ are the homogeneous components of degree $m$ of $P_{\pi}$, $P_{\tilde\pi'}$, respectively, and $\nu$ is the socle degree of each of $A$ and $\tilde A$ (observe that the socle degrees of $A$ and $\tilde A$ are equal since by (\ref{lineq1}) one has $\deg P_{\pi}=\deg P_{\tilde\pi'}$). 

Let $\pi_m$ be the symmetric $\Soc(A)$-valued $m$-form on ${\mathcal K}$ defined as follows:
\begin{equation}
\pi_m(x_1,\dots,x_m):=\pi(x_1\cdots x_m),\quad x_1,\dots,x_m\in{\mathcal K},\,\,m=2,\dots,\nu.\label{defforms}
\end{equation} 
By (\ref{poly}) we have 
\begin{equation}
P_{\pi}^{[m]}(x)=\frac{1}{m!}\pi_m(x,\dots,x),\quad x\in{\mathcal K},\,\, m=2,\dots,\nu.\label{connect1}
\end{equation}
We will focus on the forms $\pi_2$ and $\pi_3$ (in fact, it is shown in Proposition 2.8 in \cite{FIKK} that every $\pi_m$ with $m>3$ is completely determined by $\pi_2$, $\pi_3$). As in \cite{FIKK}, we define a commutative product $(x,y)\mapsto x*y$ on ${\mathcal K}$ by requiring the identity
\begin{equation}
\pi_2(x*y,z)=\pi_3(x,y,z)\label{defprod}
\end{equation}
to hold for all $x,y,z\in{\mathcal K}$. Owing to the non-degeneracy of the form $b_{\pi}$ defined in (\ref{bpi}), the form $\pi_2$ is non-degenerate and therefore for any $x,y\in{\mathcal K}$ the element $x*y\in{\mathcal K}$ is uniquely determined by (\ref{defprod}). 

We need the following lemma.

\begin{lemma}\label{newprod}\sl For any two elements $x+u$, $y+v$ of ${\mathfrak m}$, with $x,y\in{\mathcal K}$ and $u,v\in\Soc(A)$, one has
\begin{equation}
(x+u)(y+v)=x*y+\pi_2(x,y).\label{sameproduct}
\end{equation}
\end{lemma}

\noindent {\bf Proof:} We write the product $xy$ with respect to the decomposition\linebreak ${\mathfrak m}={\mathcal K}\oplus\Soc(A)$ as $xy=(xy)_1+(xy)_2$, where $(xy)_1\in{\mathcal K}$ and $(xy)_2\in\Soc(A)$. It is then clear that $(xy)_2=\pi_2(x,y)$. Therefore, to prove (\ref{sameproduct}) we need to show that $x*y=(xy)_1$. This identity is obtained by a simple calculation similar to the one that occurs in the proof of Proposition 2.8 of \cite{FIKK} (cf.~Proposition 5.13 in \cite{FK}). Indeed, from (\ref{defforms}), (\ref{defprod}) for any $z\in{\mathcal K}$ one has
$$
\begin{array}{l}
\pi_2(x*y,z)=\pi_3(x,y,z)=\pi(xyz)=\\
\vspace{-0.1cm}\\
\hspace{4cm}\pi\Bigl([(xy)_1+(xy)_2]z\Bigr)=\pi((xy)_1z)=\pi_2((xy)_1,z).
\end{array}
$$
Since $\pi_2$ is non-degenerate, the above identity implies $x*y=(xy)_1$ as required.\qed
\vspace{0.1cm}\\

We shall now complete the proof of Theorem \ref{equivalence}. Let $\tilde\pi_m$ be the symmetric $\Soc(\tilde A)$-valued forms on $\tilde{\mathcal K}$ arising from the admissible projection $\tilde\pi'$ on $\tilde{\mathfrak m}$ as in (\ref{defforms}), with $m=2,\dots,\nu$. By (\ref{lineq}), (\ref{connect1}) we have
\begin{equation}
\begin{array}{l}
L_2(\pi_m(x_1,\dots,x_m))=\tilde\pi_m(L_1(x_1),\dots, L_1(x_m)),\\ 
\vspace{-0.1cm}\\
\hspace{6cm}x_1,\dots,x_m\in{\mathcal K},\,\,m=2,\dots,\nu.
\end{array}\label{impidentities}
\end{equation}
Denote by $\tilde *$ the product on $\tilde{\mathcal K}$ defined by $\tilde\pi_2$, $\tilde\pi_3$ as in (\ref{defprod}). Now, for all $x,y\in{\mathcal K}$ and $u,v\in\Soc(A)$, Lemma \ref{newprod} and identity (\ref{impidentities}) with $m=2$ yield
\begin{equation}
\makebox[250pt]{$\begin{array}{l}
{\mathcal L}\Bigl((x+u)(y+v)\Bigr)={\mathcal L}(x*y+\pi_2(x,y))=L_1(x*y)+L_2(\pi_2(x,y)),\\
\vspace{0cm}\\
{\mathcal L}(x+u){\mathcal L}(y+v)=(L_1(x)+L_2(u))(L_1(y)+L_2(v))=\\
\vspace{-0.3cm}\\
\hspace{2cm}L_1(x)\tilde *L_1(y)+\tilde\pi_2(L_1(x),L_1(y))=L_1(x)\tilde *L_1(y)+L_2(\pi_2(x,y)).
\end{array}$}\label{imporeq}
\end{equation}
Next, for any $z\in{\mathcal K}$, from (\ref{defprod}) and identity (\ref{impidentities}) with $m=2,3$ one obtains
$$
\begin{array}{l}
\tilde\pi_2(L_1(x)\tilde*L_1(y),L_1(z))=\tilde\pi_3(L_1(x),L_1(y),L_1(z))=\\
\vspace{-0.1cm}\\
\hspace{3cm}L_2(\pi_3(x,y,z))=L_2(\pi_2(x*y,z))=\tilde\pi_2(L_1(x*y),L_1(z)),
\end{array}
$$
which implies $L_1(x)\tilde*L_1(y)=L_1(x*y)$. It then follows from (\ref{imporeq}) that ${\mathcal L}\Bigl((x+u)(y+v)\Bigr)={\mathcal L}(x+u){\mathcal L}(y+v)$, that is, ${\mathcal L}:{\mathfrak m}\ra\tilde{\mathfrak m}$ is an algebra isomorphism. 

Finally, we will obtain the last statement of the theorem. Let ${\mathcal A}:{\mathfrak m}\ra\tilde{\mathfrak m}$ be a linear isomorphism such that ${\mathcal A}(S_{\pi})=S_{\tilde\pi}$. Consider the linear automorphisms ${\mathcal F}$ of ${\mathfrak m}$ and $\tilde{\mathcal F}$ of $\tilde{\mathfrak m}$ defined by
$$
\begin{array}{ll}
{\mathcal F}(x+u):=x-u,&x\in{\mathcal K},\,u\in\Soc(A),\\
\vspace{-0.3cm}\\
\tilde{\mathcal F}(x+u):=x-u,&x\in\ker\tilde\pi\cap\tilde{\mathfrak m},\,u\in\Soc(\tilde A).
\end{array}
$$
Then the composition ${\mathcal A}_0:=\tilde{\mathcal F}\circ{\mathcal A}\circ{\mathcal F}$ is a linear transformation from ${\mathfrak m}$ to $\tilde{\mathfrak m}$ that maps ${\mathcal S}_{\pi}$ onto ${\mathcal S}_{\tilde\pi}$. By the above argument, it follows that ${\mathcal A}_0$ is an algebra isomorphism. On the other hand, one has ${\mathcal A}_0={\mathcal A}$. Thus, ${\mathcal A}$ is an algebra isomorphism, and the proof of the theorem is complete.\qed

\section{Example of application of Theorem \ref{equivalence}}\label{app}
\setcounter{equation}{0}

Theorem \ref{equivalence} is particularly useful when at least one of the hypersurfaces $S_{\pi}$, $S_{\tilde\pi}$ is affinely homogeneous (recall that a subset ${\mathcal S}$ of a vector space $V$ is called affinely homogeneous if for every pair of points $p,q\in{\mathcal S}$ there exists a bijective affine map ${\mathcal A}$ of $V$ such that ${\mathcal A}({\mathcal S})={\mathcal S}$ and ${\mathcal A}(p)=q$). In this case the hypersurfaces $S_{\pi}$, $S_{\tilde\pi}$ are affinely equivalent if and only if they are linearly equivalent. Indeed, if, for instance, $S_{\pi}$ is affinely homogeneous and ${\mathcal A}:{\mathfrak m}\ra\tilde{\mathfrak m}$ is an affine equivalence between $S_{\pi}$, $S_{\tilde\pi}$, then ${\mathcal A}\circ {\mathcal A}'$ is a linear equivalence between $S_{\pi}$, $S_{\tilde\pi}$, where ${\mathcal A}'$ is an affine automorphism of $S_{\pi}$ such that ${\mathcal A}'(0)={\mathcal A}^{-1}(0)$. Clearly, in this case $S_{\tilde\pi}$ is affinely homogeneous as well. 

The proof of Theorem \ref{equivalence} shows that every linear equivalence ${\mathcal L}$ between $S_{\pi}$, $S_{\tilde\pi}$ has the block-diagonal form with respect to the decompositions\linebreak ${\mathfrak m}={\mathcal K}\oplus\Soc(A)$ and $\tilde{\mathfrak m}=\tilde{\mathcal K}\oplus\Soc(\tilde A)$, where $\tilde{\mathcal K}:=\ker\tilde\pi\cap\tilde{\mathfrak m}$, that is, there exist linear isomorphisms $L_1: {\mathcal K}\ra\tilde{\mathcal K}$ and $L_2:\Soc(A)\ra\Soc(\tilde A)$ such that ${\mathcal L}(x+u)=L_1(x)+L_2(u)$, with $x\in{\mathcal K}$, $u\in\Soc(A)$. Therefore, analogously to (\ref{lineq}), for the corresponding polynomial maps $P_{\pi}$ and $P_{\tilde\pi}$ (see (\ref{poly})) we have
\begin{equation}
L_2\circ P_{\pi}^{[m]}=P_{\tilde\pi}^{[m]}\circ L_1\quad\hbox{for all $m\ge 2$},\label{lineq2}
\end{equation}
where, as before, $P_{\pi}^{[m]}$, $P_{\tilde\pi}^{[m]}$ are the homogeneous components of degree $m$ of $P_{\pi}$, $P_{\tilde\pi}$, respectively.

Thus, Theorem \ref{equivalence} yields the following corollary (cf.~Theorem 2.11 in \cite{FIKK}) .

\begin{corollary}\label{cor}\sl Let $A$, $\tilde A$ be Gorenstein algebras of finite vector space dimension greater than 1 and socle degree $\nu$ over an infinite field $k$ with either $\char(k)=0$ or $\char(k)>\nu$. Let, further, $\pi$, $\tilde\pi$ be admissible projections on $A$, $\tilde A$, respectively.

\noindent {\rm (i)} If $A$ and $\tilde A$ are isomorphic and at least one of $S_{\pi}$, $S_{\tilde\pi}$ is affinely homogeneous, then for some linear isomorphisms $L_1: {\mathcal K}\ra\tilde{\mathcal K}$ and $L_2:\Soc(A)\ra\Soc(\tilde A)$ identity {\rm (\ref{lineq2})} holds. In this case both $S_{\pi}$ and $S_{\tilde\pi}$ are affinely homogeneous.

\noindent {\rm (ii)} If for some linear isomorphisms $L_1: {\mathcal K}\ra\tilde{\mathcal K}$ and $L_2:\Soc(A)\ra\Soc(\tilde A)$ identity {\rm (\ref{lineq2})} holds, then the hypersurfaces $S_{\pi}$, $S_{\tilde\pi}$ are linearly equivalent and therefore the algebras $A$ and $\tilde A$ are isomorphic. 
\end{corollary}

As shown in Section 8.2 in \cite{FK}, the hypersurface $S_{\pi}$ need not be affinely homogeneous in general. In \cite{Is} (see also \cite{FK}) we found a criterion for the affine homogeneity of some (hence every) hypersurface $S_{\pi}$ arising from an Artinian Gorenstein algebra $A$. Namely, $S_{\pi}$ is affinely homogeneous if and only if the action of the automorphism group of the nilpotent algebra ${\mathfrak m}$ on the set of all hyperplanes in ${\mathfrak m}$ complementary to $\Soc(A)$ is transitive. Furthermore, we showed that this condition is satisfied if $A$ is non-negatively graded in the sense that it can be represented as a direct sum 
\begin{equation}
A=\bigoplus_{j\ge0}A^j,\quad A^{j}A^{\ell}\subset A^{j+\ell},\label{sum}
\end{equation}
where $A^{j}$ are linear subspaces of $A$, with $A^0=k$ (in this case ${\mathfrak m}=\oplus_{j>0}A^j$ and $\Soc(A)=A^d$ for $d:=\max\{j:A^{j}\ne 0\}$). We stress that the grading $\{A^j\}$ in the above statement is {\it not}\, required to be standard, i.e.~$A^j$ may not coincide with $(A^1)^j$.

Note, however, that the existence of a non-negative grading on $A$ is not a necessary condition for the affine homogeneity of $S_{\pi}$. For example, for any Artinian Gorenstein algebra of socle degree not exceeding 4 the hypersurfaces $S_{\pi}$ are affinely homogeneous (see Proposition 6.5 in \cite{FK}), whereas not every such algebra admits a non-negative grading. The case $\nu=3$ is relatively easy; in fact, all algebras with this property were completely described in Theorem 4.1 of \cite{ER}. The case $\nu=4$ is much harder (see Section 7.2 in \cite{FK}), and the affine homogeneity of $S_{\pi}$ makes Corollary \ref{cor} an important tool in this case. We are confident that the approach discussed in this paper will help make significant advances on the classification problem at least for $\nu=4$. Note that, as shown in Proposition 7.5 of \cite{FK}, the classification result of Theorem 4.1 of \cite{ER} can indeed be obtained by utilizing this method.

We will now give an example of how one can hope to apply our technique to algebras of socle degree 4.

\begin{example}\label{appl}\rm Let $A_0$ be the Artinian Gorenstein algebra introduced at the end of Section 8.1 in \cite{FK}, namely,
$$
A_0:=k[[x,y,z]]/J(x^4+xy^2+y^3+xz^2),
$$
where $k[[x,y,z]]$ is the algebra of formal power series in $x,y,z$ and $J(f)$ is the Jacobian ideal of $f$, i.e.~the ideal generated by the first-order partial derivatives of $f$. As explained in \cite{FK}, this algebra does not admit {\it any}\, non-negative grading. Furthermore, the Hilbert function of $A_0$ is $\{1,3,3,1,1\}$, hence the associated graded algebra of $A_0$ is not Gorenstein (cf.~Proposition 9 in \cite{W}). Also, one has $\nu=4$ and $\dim_kA_0=9$. 

Consider the following monomials in $k[[x,y,z]]$:
\begin{equation}
x^4,\,\, x,\,\, x^2,\,\, x^3,\,\, y,\,\, z,\,\, yz,\,\, z^2\label{monomials}
\end{equation}
and let $e_{0,0},\dots,e_{0,7}$ be the vectors in the maximal ideal ${\mathfrak m}_0$ of $A_0$ represented, respectively, by these monomials. It is not hard to show that $e_{0,0},\dots,e_{0,7}$ are linearly independent, hence they form a basis of ${\mathfrak m}_0$. Choose $\pi_0$ to be the admissible projection on $A_0$ defined by the condition 
$$
\ker\pi_0\,\cap\,{\mathfrak m}_0=\langle e_{0,1},\dots,e_{0,7}\rangle=:{\mathcal K}_0,
$$
where $\langle{\,}\cdot{\,}\rangle$ denotes linear span (observe that $\Soc(A_0)=\langle e_{0,0}\rangle$). Then, letting $x_1,\dots,x_7$ be the coordinates in ${\mathcal K}_0$ with respect to the basis $e_{0,1},\dots,e_{0,7}$ and identifying $\Soc(A_0)$ with $k$ by means of $e_{0,0}$, we compute:\footnote{Here and below the relevant polynomials $P_{\pi}$ were found by using a {\tt Singular}-based computer program, which had been kindly made available to us by W. Kaup.}
\begin{equation}
\begin{array}{l}
\displaystyle P_0:=P_{\pi_{{}_0}}=x_1x_3+\frac{1}{2}x_2^2+6x_2x_4-\frac{8}{3}x_4x_7-\frac{8}{3}x_5x_6+\\
\vspace{-0.3cm}\\
\hspace{3cm}\displaystyle\frac{1}{2}x_1^2x_2+3x_1^2x_4-2x_1x_4^2+\frac{4}{9}x_4^3-\frac{4}{3}x_4x_5^2+\frac{1}{24}x_1^4.
\end{array}\label{p0}
\end{equation}

We will now perturb the algebra $A_0$ as follows:
$$
A_t:=k[[x,y,z]]/J(x^4+tx^5+xy^2+y^3+xz^2),\quad t\in k^*.
$$
It is not hard to check that for every $t$ the algebra $A_t$ is Artinian Gorenstein of vector space dimension 9 and socle degree 4. By utilizing part (ii) of Corollary \ref{cor}, we will now show that $A_t$ is in fact isomorphic to $A_0$ for all $t$.

Let $e_{t,0},\dots,e_{t,7}$ be the basis in the maximal ideal ${\mathfrak m}_t$ of $A_t$ whose elements are represented, respectively,  by monomials (\ref{monomials}). Let $\pi_t$ be the admissible projection on $A_t$ defined by
$$
\ker\pi_t\,\cap\,{\mathfrak m}_t=\langle e_{t,1},\dots,e_{t,7}\rangle=:{\mathcal K}_t
$$
(observe that $\Soc(A_t)=\langle e_{t,0}\rangle$). Then, denoting by $y_1,\dots,y_7$ the coordinates in ${\mathcal K}_t$ with respect to the basis $e_{t,1},\dots,e_{t,7}$ and identifying $\Soc(A_t)$ with $k$ by means of $e_{t,0}$, we have
\begin{equation}
\begin{array}{l}
\hspace{-0.5cm}\displaystyle P_t:=P_{\pi_{{}_t}}=y_1y_3+\frac{1}{2}y_2^2+6y_2y_4-\frac{8}{3}y_4y_7-\frac{8}{3}y_5y_6+\frac{15t}{2}y_1y_4-\frac{5t}{2}y_4^2+\\
\vspace{-0.3cm}\\
\hspace{3.5cm}\displaystyle\frac{1}{2}y_1^2y_2+3y_1^2y_4-2y_1y_4^2+\frac{4}{9}y_4^3-\frac{4}{3}y_4y_5^2+\frac{1}{24}y_1^4.
\end{array}\label{pt}
\end{equation}
Thus, we see that in the coordinates chosen as above the polynomials $P_0$ and $P_t$ coincide in homogeneous components of degrees 3 and 4 but their quadratic terms differ.

We now identify $\Soc(A_0)$ and $\Soc( A_t)$ by means of the vectors $e_{0,0}$ and $e_{t,0}$. Then, by part (ii) of Corollary \ref{cor}, to prove that $A_t$ is isomorphic to $A_0$, it suffices to find a linear isomorphism $L: {\mathcal K}_t\ra{\mathcal K}_0$ such that 
\begin{equation}
P_t^{[m]}=P_0^{[m]}\circ L\quad\hbox{for $m=2,3,4$},\label{mainex}
\end{equation}
where, as before, $P_0^{[m]}$, $P_t^{[m]}$ are the homogeneous components of degree $m$ of $P_0$, $P_t$, respectively. Define $L$ by 
$$
\begin{array}{l}
x_j=y_j,\quad\hbox{for $j\ne 3,7$},\\
\vspace{-0.3cm}\\
\displaystyle x_3=y_3+\frac{15t}{2}y_4, \\
\vspace{-0.3cm}\\
\displaystyle x_7=\frac{15t}{16}y_4+y_7. 
\end{array}
$$
Clearly, this linear transformation satisfies (\ref{mainex}), which shows that $A_t$ is indeed isomorphic to $A_0$ for every $t$ as claimed.

By the last statement of Theorem \ref{equivalence}, the map $L$ together with the identification of $\Soc(A_0)$ and $\Soc( A_t)$ is in fact an algebra isomorphism between ${\mathfrak m}_t$ and ${\mathfrak m}_0$. More precisely, the map defined on the basis elements as
$$
\begin{array}{l}
e_{t,j}\mapsto e_{0,j},\quad \hbox{for $j\ne 4$},\\
\vspace{-0.3cm}\\
\displaystyle e_{t,4}\mapsto \frac{15t}{2}e_{0,3}+e_{0,4}+\frac{15t}{16}e_{0,7}
\end{array}
$$
is an isomorphism from ${\mathfrak m}_t$ onto ${\mathfrak m}_0$. The corresponding isomorphism between $A_t$ and $A_0$ is induced by the automorphism of $k[[x,y,z]]$ given by the following change of variables:
\begin{equation}
x\mapsto x,\quad y\mapsto y+\frac{15t}{16}z^2+\frac{15t}{2}x^3,\quad z\mapsto z.\label{powerseraut}
\end{equation} 
\end{example}
It would be interesting to see whether formula (\ref{powerseraut}) could be obtained directly using ideal generators. 

One can produce an alternative proof of isomorphism of $A_t$ and $A_0$ by making use of Macaulay inverse systems. We will explore this possibility in the next section and compare it with the proof given above.

\section{Nil-polynomials and Macaulay\\ inverse systems}\label{invsys}
\setcounter{equation}{0}

With the exception of a further discussion of Example \ref{appl} given below, the material of this section is contained in \cite{AI}. Since this material is highly relevant to the present paper's theme and the proofs involved are not long, we reproduce it here for the reader's benefit.

As before, let $A$ be an Artinian Gorenstein algebra over a field $k$, with maximal ideal ${\mathfrak m}$ and socle degree $\nu$, where we assume that $\dim_{k}A>1$ and either $\char(k)=0$ or $\char(k)>\nu$. Next, let 
$$
M:=\emb\dim A:=\dim_{k}{\mathfrak m}/{\mathfrak m}^2
$$
be the embedding dimension of $A$ (notice that $M\ge 1$). Choose a basis ${\mathcal B}=\{e_1,\dots,e_M\}$ in a complement to ${\mathfrak m}^2$ in ${\mathfrak m}$ and fix a linear form\linebreak $\omega:A\to k$ with kernel complementary to $\Soc(A)$. Now, we introduce the following polynomial:
\begin{equation}
Q_{\omega,{\mathcal B}}(x_1,\dots,x_M):=\sum_{j=0}^{\nu}\frac{1}{j!}\omega\Bigl((x_1e_1+\dots+x_Me_M)^j\Bigr).\label{restrnilpol}
\end{equation} 
Notice that in (\ref{restrnilpol}) the element $(x_1e_1+\dots+x_Me_M)^j\in A$ may have a non-trivial projection to $\Soc(A)$ parallel to $\ker\omega$ even for $j<\nu$, in which case $\omega(x_1e_1+\dots+x_Me_M)^j\ne 0$. In formulas (\ref{q0}) and (\ref{qt}) below we compute polynomials of this kind for the algebras $A_0$ and $A_t$ from Example \ref{appl}.

Further, the elements $e_1,\dots,e_M$ generate $A$ as an algebra, hence $A$ is isomorphic to $k[x_1,\dots,x_M]/I$, where $I$ is the ideal of all relations among $e_1,\dots,e_M$, i.e.~polynomials $f\in k[x_1,\dots,x_M]$ with $f(e_1,\dots,e_M)=0$. Observe that $I$ contains the monomials $x_1^{\nu+1},\dots,x_M^{\nu+1}$. 

From now on we assume that $\char(k)=0$ (cf.~Remark \ref{fieldsofposchar} below). For $f,g\in k[x_1,\dots,x_M]$ define
$$
 f\star g:=f\left(\frac{\partial}{\partial x_1},\dots,\frac{\partial}{\partial x_M}\right)(g).
$$
It is well-known that, since the quotient $k[x_1,\dots,x_M]/I$ is Gorenstein, there is a polynomial $g\in k[x_1,\dots,x_M]$ of degree $\nu$ such that $I=\Ann(g)$, where
$$
\Ann(g):=\left\{f\in k[x_1,\dots,x_M]: f\star g=0\right\}
$$
is the {\it annihilator}\, of $g$. The freedom in choosing $g$ with $\Ann(g)=I$ is fully understood, namely, if $g_1,g_2\in k[x_1,\dots,x_M]$ satisfy $\Ann(g_1)=\Ann(g_2)=I$, then $g_1=h\star g_2$, where $h\in k[x_1,\dots,x_M]$ with $h(0)\ne 0$. Any polynomial $g$ with  $I=\Ann(g)$ is called a {\it Macaulay inverse system}\, for the Artinian Gorenstein quotient $k[x_1,\dots,x_M]/I$.

Conversely, for any non-zero element $g\in k[x_1,\dots,x_M]$ the quotient $k[x_1,\dots,x_M]/\Ann(g)$ is a (local) Artinian Gorenstein algebra of socle degree $\deg g$, hence any polynomial is an inverse system of some Artinian Gorenstein quotient. It is well-known that inverse systems can be used for solving the isomorphism problem for quotients of this kind as explained, for example, in Proposition 2.2 in \cite{ER} (a precise statement is given at the end of this section). For details on inverse systems we refer the reader to \cite{M}, \cite{Em}, \cite{Ia} (a brief survey given in \cite{ER} is also helpful).

In applications it is desirable to have an explicit formula for computing an inverse system for any Artinian Gorenstein quotient. As the following theorem shows, formula (\ref{restrnilpol}) achieves this purpose.

\begin{theorem}\label{main}{\rm\cite{AI}}\,\sl Let $A$ be a Gorenstein algebra of finite vector space dimension greater than 1 over a field of characteristic zero, with maximal ideal ${\mathfrak m}$, socle degree $\nu$ and embedding dimension $M$. Choose a basis\linebreak ${\mathcal B}=\{e_1,\dots,e_M\}$ in a complement to ${\mathfrak m}^2$ in ${\mathfrak m}$, fix a linear form $\omega:A\to k$ with kernel complementary to $\Soc(A)$ and consider the polynomial $Q_{\omega,{\mathcal B}}$ as defined in {\rm (\ref{restrnilpol})}. Further, using the basis ${\mathcal B}$ write the algebra $A$ as a quotient $k[x_1,\dots,x_M]/I$. Then $Q_{\omega,{\mathcal B}}$ is an inverse system for $k[x_1,\dots,x_M]/I$.
\end{theorem}

\noindent{\bf Proof:} Fix any polynomial $f\in k[x_1,\dots,x_M]$
$$
f=\sum_{0\le i_1,\dots,i_M\le N}a_{i_1,\dots,i_M}x_1^{i_1}\dots x_M^{i_M}
$$ 
and calculate
\begin{equation}
\begin{array}{l}
\displaystyle  f\star Q_{\omega,{\mathcal B}}=f\left(\frac{\partial}{\partial x_1},\dots,\frac{\partial}{\partial x_M}\right)(Q_{\omega,{\mathcal B}})=\\
\vspace{-0.3cm}\\
\displaystyle\sum_{0\le i_1,\dots,i_M\le N}a_{i_1,\dots,i_M}\sum_{j= i_1+\dots+i_M}^{\nu}\frac{1}{(j-(i_1+\dots+i_M))!}\times\\
\vspace{-0.45cm}\\
\hspace{2.2cm}\omega\Bigl((x_1e_1+\dots+x_Me_M)^{j-(i_1+\dots+i_M)}e_1^{i_1}\dots e_M^{i_M}\Bigr)=\\
\vspace{-0.3cm}\\
\displaystyle\sum_{m=0}^{\nu}\frac{1}{m!}\omega\Bigl((x_1e_1+\dots+x_Me_M)^m\times\\
\vspace{-0.7cm}\\
\hspace{4cm}\displaystyle\sum_{\mbox{\tiny$\begin{array}{l} 0\le i_1,\dots,i_M\le N,\\\vspace{-0.15cm}\\ i_1+\dots+i_M\le\nu-m\end{array}$}}
a_{i_1,\dots,i_M}e_1^{i_1}\dots e_M^{i_M}\Bigr)=\\
\vspace{-0.3cm}\\
\displaystyle\sum_{m=0}^{\nu}\frac{1}{m!}\omega\Bigl((x_1e_1+\dots+x_Me_M)^m\,f(e_1,\dots,e_M)\Bigr).
\end{array}\label{diff}
\end{equation}
Since $I$ is the ideal of all relations among $e_1,\dots,e_M$, formula (\ref{diff}) implies $I\subset\Ann(Q_{\omega,{\mathcal B}})$.

Conversely, let $f\in k[x_1,\dots,x_M]$ be an element of $\Ann(Q_{\omega,{\mathcal B}})$. Then (\ref{diff}) yields
\begin{equation}
\sum_{m=0}^{\nu}\frac{1}{m!}\omega\Bigl((x_1e_1+\dots+x_Me_M)^m\,f(e_1,\dots,e_M)\Bigr)=0.\label{eq2}
\end{equation}
Collecting the terms containing $x_1^{i_1}\dots x_M^{i_M}$ in (\ref{eq2}) we obtain
\begin{equation}
\omega\Bigl(e_1^{i_1}\dots e_M^{i_M}\,f(e_1,\dots,e_M)\Bigr)=0\label{eq3}
\end{equation}
for all indices $i_1,\dots,i_M$. Since $e_1,\dots,e_M$ generate $A$, identities (\ref{eq3}) yield 
\begin{equation}
\omega\Bigl(A\, f(e_1,\dots,e_M)\Bigr)=0.\label{nondeg}
\end{equation}
Further, since the bilinear form $(a,b)\mapsto \omega(ab)$ is non-degenerate on $A$, identity (\ref{nondeg}) implies $f(e_1,\dots,e_M)=0$. Therefore $f\in I$, which shows that $I=\Ann(Q_{\omega,{\mathcal B}})$ as required.\qed
\vspace{0.1cm}\\

We will now make a number of useful remarks.

\begin{remark}\label{fieldsofposchar}\rm For simplicity, we have chosen to discuss inverse systems only under the assumption $\char(k)=0$. For fields of positive characteristics one needs to pass to divided power rings (see, e.g.~Lemma 1.2 in \cite{Ia}). If\linebreak $\char(k)>\nu$, one can naturally think of $Q_{\omega,{\mathcal B}}$ as an element of the corresponding divided power ring, and Theorem \ref{main} remains correct.
\end{remark}

\begin{remark}\label{nilpolinvsys}\rm Fix a hyperplane $\Pi$ in ${\mathfrak m}$ complementary to $\Soc(A)$. A $k$-valued polynomial $P$ on $\Pi$ is called a {\it nil-polynomial for $A$}\, if there exists a linear form $\rho:A\ra k$ such that $\ker\rho=\langle\Pi,{\bf 1}\rangle$ and $P=\rho\circ\exp|_{\Pi}$. Note that $\deg P=\nu$. Upon identification of $\Soc(A)$ with $k$, the class of nil-polynomials coincides with that of $\Soc(A)$-valued polynomial maps $P_{\pi}$ introduced in (\ref{poly}). If $\dim_{k}A>2$, then $\nu\ge 2$ and $\Pi$ contains an $M$-dimensional subspace that forms a complement to ${\mathfrak m}^2$ in ${\mathfrak m}$. Fix any such subspace $V$, choose a basis ${\mathcal B}=\{e_1,\dots,e_M\}$ in $V$ and let $x_1,\dots,x_M$ be the coordinates in $V$ with respect to this basis. Then the polynomial $Q_{\rho,{\mathcal B}}$ given by formula (\ref{restrnilpol}) is exactly the restriction of the nil-polynomial $P$ to $V$ written in the coordinates $x_1,\dots,x_M$. Thus, one way to explicitly obtain an inverse system for an Artinian Gorenstein quotient is to restrict a socle-valued map $P_{\pi}$ to a complement to ${\mathfrak m}^2$ in ${\mathfrak m}$ lying in $\ker\pi\cap{\mathfrak m}$ and identify the socle with the field $k$. This observation provides a link between the classical approach to Artinian Gorenstein algebras by means of inverse systems and our criterion in Theorem \ref{equivalence}.   
\end{remark}

\begin{remark}\label{stgradedhom}\rm Suppose that $A$ is a {\it standard graded algebra}, i.e.~$A$ can be represented in the form (\ref{sum}) with $A^j=(A^1)^j$ for $j\ge 1$. Set
$$
\Pi:=\bigoplus_{j=1}^{\nu-1}A^j,\quad V:=A^1.
$$
In this case, for any choice of a basis ${\mathcal B}=\{e_1,\dots,e_M\}$ in $V$, the ideal $I$ is homogeneous, i.e.~generated by homogeneous relations. For an arbitrary nil-polynomial $P$ on $\Pi$ its restriction $Q_{\rho,{\mathcal B}}$ to $V$ coincides with the homogeneous component of degree $\nu$ of $P$:
$$
Q_{\rho,{\mathcal B}}(x_1,\dots,x_M)=\frac{1}{\nu!}\omega\Bigl((x_1e_1+\dots+x_Me_M)^{\nu}\Bigr).
$$
Thus, Theorem \ref{main} yields a simple proof of the well-known fact that a standard graded Artinian Gorenstein algebra, when written as a quotient by a homogeneous ideal, admits a homogeneous inverse system and all homogeneous inverse systems for such algebras are mutually proportional (see \cite{Em} and  pp.~79--80 in \cite{M}).
\end{remark}

\begin{remark}\label{remembdim}\rm Theorem \ref{main} easily generalizes to the case of Artinian Gorenstein quotients $k[x_1,\dots,x_m]/I$, where $I$ lies in the ideal generated by\linebreak $x_1,\dots,x_m$ and $m$ is not necessarily equal to the embedding dimension $M$ of the quotient. Indeed, let $e_1,\dots,e_m$ be the elements of $k[x_1,\dots,x_m]/I$ represented by $x_1,\dots,x_m$, respectively, and consider the element of $k[x_1,\dots,x_m]$ defined as follows:
\begin{equation}
R(x_1,\dots,x_m):=\sum_{j=0}^{\nu}\frac{1}{j!}\omega\Bigl((x_1e_1+\dots+x_me_m)^j\Bigr),\label{explinvss}
\end{equation}
where $\omega$ is a linear form on $k[x_1,\dots,x_m]/I$ with kernel complementary to the socle and $\nu$ is the socle degree of $k[x_1,\dots,x_m]/I$. Then, arguing as in the proof of Theorem \ref{main}, we see that $R$ is an inverse system for $k[x_1,\dots,x_m]/I$. Thus, (\ref{explinvss}) is an explicit formula providing an inverse system for any Artinian Gorenstein quotient, and any other inverse system is obtained as $h\star R$, where $h\in k[x_1,\dots,x_m]$ does not vanish at the origin. Notice that for\linebreak $m>M$, no inverse system as in (\ref{explinvss}) comes from restricting a nil-polynomial to a subspace of ${\mathfrak m}$ complementary to ${\mathfrak m}^2$.
\end{remark}

Thus, in order to decide whether two Artinian Gorenstein algebras are isomorphic, one can use either the classical approach, which utilizes inverse systems, or our method, which relies on nil-polynomials, and the two techniques are related as explained in Theorem \ref{main} and Remark \ref{nilpolinvsys}. In a particular situation one of the approaches may work better than the other. For instance, as we saw in the preceding section, the technique based on nil-polynomials is very appropriate for establishing that the algebras $A_t$ in Example \ref{appl} are all isomorphic to $A_0$. We will now obtain a different proof of this statement by the method based on inverse systems.

We first describe this method in general following the discussion that precedes Proposition 2.2 in \cite{ER}. For $j=1,2$, let $B_j:=k[x_1,\dots,x_m]/I_j$ be an Artinian Gorenstein quotient with socle degree $\nu\ge 1$ and inverse system $g_j$. Denote by $k_{\nu}[x_1,\dots,x_m]$ the vector space of polynomials in $k[x_1,\dots,x_m]$ of degree not exceeding $\nu$. If $T_1,\dots,T_m\in k_{\nu}[x_1,\dots,x_m]$ vanish at the origin and have linearly independent linear parts, they induce a linear automorphism $\Phi_{T_1,\dots,T_m}$ of $k_{\nu}[x_1,\dots,x_m]$ as follows:
$$
\Phi_{T_1,\dots,T_m}: f(x)\mapsto f(T_1(x),\dots,T_m(x))\,\,\hbox{(mod terms of degree $>\nu$)},
$$
where $f\in k_{\nu}[x_1,\dots,x_m]$ and $x:=(x_1,\dots,x_m)$. We now introduce a symmetric $k$-valued bilinear form on $k_{\nu}[x_1,\dots,x_m]$ as
$$
[f,\tilde f]:=(f\star \tilde f)(0),\quad f,\tilde f\in k_{\nu}[x_1,\dots,x_m].
$$
This form is clearly non-degenerate, and one can consider the linear map adjoint to $\Phi_{T_1,\dots,T_m}$, i.e.~the automorphism $\Phi_{T_1,\dots,T_m}^*$ of $k_{\nu}[x_1,\dots,x_m]$ defined by requiring that the identity
$$
[f,\,\,\Phi_{T_1,\dots,T_m}^*(\tilde f)]=[\Phi_{T_1,\dots,T_m}(f),\,\,\tilde f]
$$
hold for all  $f,\tilde f\in k_{\nu}[x_1,\dots,x_m]$.

Then, as explained in detail in \cite{ER}, the algebras $B_1$ and $B_2$ are isomorphic if and only if there exist $T_1,\dots,T_m$ as above and $h\in k[x_1,\dots,x_m]$ with $h(0)\ne 0$ such that
\begin{equation}
\Phi_{T_1,\dots,T_m}^*(g_1)=h\star g_2.\label{invsysequiv}
\end{equation}
Note that the right-hand side of formula (\ref{invsysequiv}) is simply a replacement of the inverse system $g_2$ by another inverse system for the algebra $B_2$. Clearly, (\ref{invsysequiv}) is equivalent to the identity
\begin{equation}
\Bigl(\Phi_{T_1,\dots,T_m}(f)\star g_1\Bigr)(0)=\Bigl(f\star(h\star g_2)\Bigr)(0)\label{mainids}
\end{equation}
being satisfied for all $f\in k_{\nu}[x_1,\dots,x_m]$.

We now return to Example \ref{appl} and consider the algebras $A_0$, $A_t$ of socle degree $\nu=4$ introduced there, with $t\in k^*$. Let $V_0:=\langle e_{0,1}, e_{0,4}, e_{0,5}\rangle$. This subspace forms a complement to ${\mathfrak m}_0^2$ in ${\mathfrak m}_0$. Then, by Remark \ref{nilpolinvsys}, setting $x_2=x_3=x_6=x_7=0$ and replacing, respectively, $x_1$, $x_4$, $x_5$ by $z_1$, $z_2$, $z_3$
in formula (\ref{p0}), we obtain the inverse system
\begin{equation}
\begin{array}{l}
\displaystyle Q_0:=3z_1^2z_2-2z_1z_2^2+\frac{4}{9}z_2^3-\frac{4}{3}z_2z_3^2+\frac{1}{24}z_1^4\label{q0}
\end{array}
\end{equation}
for the algebra $A_0$ represented as the quotient $k[z_1,z_2,z_3]/I_0$, with $I_0$ being the ideal of all relations among $e_{0,1}$, $e_{0,4}$, $e_{0,5}$.

Analogously, let $V_t:=\langle e_{t,1}, e_{t,4}, e_{t,5}\rangle$. This subspace is a complement to ${\mathfrak m}_t^2$ in ${\mathfrak m}_t$. Setting $y_2=y_3=y_6=y_7=0$ and replacing, respectively, $y_1$, $y_4$, $y_5$ by $z_1$, $z_2$, $z_3$ in formula (\ref{pt}), one obtains the inverse system
\begin{equation}
\begin{array}{l}
\displaystyle Q_t:=\frac{15t}{2}z_1z_2-\frac{5t}{2}z_2^2+3z_1^2z_2-2z_1z_2^2+\frac{4}{9}z_2^3-\frac{4}{3}z_2z_3^2+\frac{1}{24}z_1^4\label{qt}
\end{array}
\end{equation}
for the algebra $A_t$ represented as the quotient $k[z_1,z_2,z_3]/I_t$, where $I_t$ is the ideal of all relations among $e_{t,1}$, $e_{t,4}$, $e_{t,5}$.

In accordance with (\ref{mainids}), to show that $A_t$ is isomorphic to $A_0$, we need to find polynomials $T_1,T_2,T_3\in k_4[z_1,z_2,z_3]$ vanishing at the origin and having linearly independent linear parts, as well as a polynomial\linebreak $h\in k[z_1,z_2,z_3]$ with $h(0)\ne 0$, such that for all $f\in k_4[z_1,z_2,z_3]$ one has
\begin{equation}
\Bigl(\Phi_{T_1,T_2,T_3}(f)\star Q_t\Bigr)(0)=\Bigl(f\star(h\star Q_0)\Bigr)(0).\label{exids}
\end{equation}
After much computational experimentation we discovered that (\ref{exids}) is satisfied for the following choice of $T_j$ and $h$:
$$
T_1=z_1,\quad T_2=z_2-\frac{15t}{16}z_3^2-\frac{15t}{2}z_1^3,\quad T_3=z_3,\quad h\equiv 1.
$$
Observe that, upon identification of $x$, $y$, $z$ with $z_1$, $z_2$, $z_3$, respectively, the change of variables $z_j\mapsto T_j$ is the inverse of the change of variables shown in (\ref{powerseraut}), which agrees with formulas (5), (6) in \cite{ER}. Thus, we have obtained a second proof of the fact that $A_t$ is indeed isomorphic to $A_0$ for all $t$.

To summarize, in the case of the algebras $A_t$ from Example \ref{appl}, the method for establishing isomorphism between Artinian Gorenstein algebras based on inverse systems requires a substantially greater computational effort than that based on nil-polynomials as presented in Section \ref{app}. Notice that nil-polynomials were utilized in this second proof as a tool for explicitly producing inverse systems.

{\obeylines
\noindent Mathematical Sciences Institute
\noindent The Australian National University
\noindent Acton, ACT 2601
\noindent Australia
\noindent e-mail: alexander.isaev@anu.edu.au
}

\end{document}